\input amstex
\input amsppt.sty
\magnification\magstep1

\def\ni\noindent
\def\sbs{\subset}

\def\diam{\operatorname{diam}}

\def\R{\text{\bf R}}

\def\Q{\text{\bf Q}}
\def\Z{\text{\bf Z}}

\def\N{\text{\bf N}}

\def\sU{\Cal U}

\hoffset= 0.0in
\voffset= 0.0in
\hsize=32pc
\vsize=38pc
\baselineskip=24pt
\NoBlackBoxes
\topmatter
\author
A.N. Dranishnikov
\endauthor

\title
Lipschitz Cohomology, Novikov conjecture, and Expanders
\endtitle

%\enddate
\abstract
We present sufficient conditions for the cohomology of a closed
aspherical manifold to be proper Lipschitz in sense of Connes-Gromov-Moscovici [CGM].
The conditions are stated in terms of the Stone-\v Cech  compactification
of the universal cover
of a manifold. We show that these conditions are formally weaker than the sufficient
conditions for the Novikov conjecture given in [CP].
Also we show that the Cayley graph of the fundamental group of a closed aspherical
manifold with proper Lipschitz cohomology cannot
contain an expander in the coarse sense. In particular, this rules
out a Lipschitz cohomology
approach to the Novikov Conjecture for
recent Gromov's examples of exotic groups.
\endabstract

\thanks The author was partially supported by NSF grant DMS-0305152.
\endthanks

\address University of Florida, Department of Mathematics, P.O.~Box~118105,
358 Little Hall, Gainesville, FL 32611-8105, USA, and Steklov Mathematical Institute,
Gubkin Str. 8, Moscow, GSP 117966, RUSSIA
\endaddress


\email  dranish\@math.ufl.edu
\endemail

\endtopmatter

\document
\head \S1 Introduction \endhead

Lipschitz cohomology classes of a group were defined by Connes, Gromov and Moscovici
[CGM]. They proved the Novikov Higher signature conjecture for such classes.
Lipschitz cohomology are defined as the images under certain slant product homomorphisms
denoted in [CGM] as $\alpha_{\cap}$. Since the analytic assembly map is a slant product
in K-theory (with the Mischenko line bundle)[FRR], it allows to make a connection between
the Lipschitz cohomology and the Novikov Conjecture (see [CGM]).
As application Connes, Gromov and Moscovici proved the Novikov conjecture for the hyperbolic groups.
In the same spirit using Lipschitz cohomology T. Kato proved the Novikov conjecture for
the combable groups [K1],[K2].
In these cases all cohomology classes of a group $\Gamma$ are Lipschitz and moreover, they
belong to the image of one homomorphism $\alpha_{\cap}$. In such situation we say that
a group $\Gamma$ has {\it canonically Lipschitz cohomology}. We consider the case when
$\Gamma$ admits a compact classifying space $B\Gamma$. In this case the cohomologies of
a group $\Gamma$ coincide with the cohomologies of the space $B\Gamma$ and hence
we can speak about
Lipschitz cohomology of a space. A typical example of a space with
the canonically Lipschitz cohomology is a closed manifold $M$ with a nonpositive sectional
curvature [CGM]. In this case the construction of corresponding homomorphism
$\alpha_{\cap}$ relies on
the following fact: Every such  manifold $M$ is aspherical, i.e.
its universal cover $X$ is
contractible, and the inverse of the exponential map $\log :X\to T_x$ is 1-Lipschitz
for every $x\in X$. In this paper we consider closed aspherical manifolds with
canonically Lipschitz cohomologies.

The main result of the paper (Theorem 1) is a reformulation of the property of
a manifold $M$ to have Lipschitz cohomology in terms of compactifications
of its universal cover $X$. Clearly that an action of $\Gamma$ on $X$ extends to an action
on the Stone-\v Cech compactification. We show that for an $n$-manifold $M=B\Gamma$ to have
canonically Lipschitz rational cohomology is equivalent to possess an essential
$\Gamma$-invariant map of the corona of the Stone-\v Cech compactification of $X$ onto the
$n-1$-dimensional sphere.
We also reformulate this condition in terms of a compactification of $X$ introduced
by Hurder [Hu]
which can be defined in simple terms (Proposition 4) as the maximal equivariant
compactification of $X$ for the fundamental group action
such that the corona is the fixed point set. Theorem 1 implies in
particular that the Novikov Conjecture holds for the group
$\pi_1(M)$ if the Hurder corona of $X$
carries a cohomology class which hits the fundamental class of $X$ under the
coboundary homomorphism $\delta:\check{H}^{n-1}(\bar X\setminus X)\to H^n_c(X)$.

There are several approaches to the Novikov Conjecture formulated in terms
of coarse geometry and coarse topology
due to Carlsson-Pedersen, Ferry-Weinberger,
Gromov, Higson-Roe, Hurder and others
(see [FRR]). It turns out that
many of them can be reduced to the Lipschitz cohomology approach
though I failed to recover Gromov's approach via
hypereuclidean manifolds ([G2], page 152) on this way.
One of the corollaries of Theorem 1 (Theorem 2) deals with a reduction
which is due to Carlsson-Pedersen.
We show that the Carlsson-Pedersen conditions on $M$ imply that $M$ has Lipschitz
cohomologies.

Perhaps the most advanced results on the Novikov Conjecture is a theorem of Gouliang Yu
which states that the Novikov Conjecture holds true for a group that admits as
a metric space with the word metric a course imbedding in a Hilbert space [Yu],[H].
Clearly, after this theorem the question by Gromov [G1] whether every metric separable space is
coarsely imbeddable in a Hilbert space became very important.
Unfortunately the answer to Gromov's question turned out to be negative.
First example of a countable
discrete metric spaces
which is not coarsely embeddable in the Hilbert space was constructed
in [DGLY]. Then Gromov noticed
that expanders form  a natural counterexample to his question.
An {\it expander} is a growing  sequence of graphs ${X_n}$ with uniformly
bounded valence and with uniform spectral gap for the combinatorial
Laplacians for the first nonzero eigenvalues:
$\lambda_1(X_n)\ge\delta>0$.

Gromov constructed finitely presented groups $\Gamma$ that  contains expander
in some coarse sense.
Moreover, Gromov constructed a closed aspherical 4-manifold $M$ with the fundamental
group containing an expander [G3],[G4]. We should note that
Gromov's imbedding conditions are weaker than those of a quasi-isometric imbedding or
a coarse imbedding. The latter means an imbedding in the coarse category [R01],[Ro2]
which is often called in the literature by a misleading name
{\it uniform imbedding} [G1] used for a different thing in the 
classical analysis. We will call here Gromov's imbedding
of expander by {\it quasi-coarse imbedding}.

Thus, expanders form an obstacle for
proving the Novikov conjecture using Yu's theorem. Also they were used to disprove
some versions of the Baum-Connes conjecture [HLS]. In this paper we show that
they are obstructions for proving the Novikov conjecture via Lipschitz cohomology
as well (Theorem 3). We show that a group with the canonically Lipschitz cohomology
cannot contain an expander.
Probably I should mention here that this fact was independently
discovered by Gromov (oral communication). Certainly, it was known in
the main partial case. Namely, Gromov [G4] and
Higson (unpublished) proved independently that the fundamental group
of a manifold with nonpositive
sectional curvature cannot contain an expander in the coarse (and Gromov's) sense.

I would like to express my thanks to Michael Gromov, Nigel
Higson, John Roe and Shmuel Weinberger for useful discussions and valuable remarks.

\head \S 2 Proper Lipschitz cohomology\endhead

In this paper all homology and cohomology will be taken with rational coefficients $\Q$.
Suppose that a discrete group $\Gamma$ acts freely and properly on a space $X$.
By $H^{\Gamma}_*(X)$ we denote the equivariant homology of $X$, i.e. the homology
of the chain complex $C^{lf}_*(X)^{\Gamma}$ of locally finite equivariant (singular) chains.
Let $q:X\to M$ be the projection onto the orbit space. Let $c$ be a (singular) simplex
of $M$ and let $\bar c$ be its lift. Then the formula $\mu(c)=\sum_{g\in\Gamma}g\bar c$
defines an isomorphism of chain complexes $\mu=\mu_q:C_*(M)\to C_*^{lf}(X)^{\Gamma}$
provided $M$ is compact. Then $\mu$ defines an isomorphism
$\mu_*:H_*(M)\to H_*^{\Gamma}(X)$. If $M$ is not compact $\mu$ defines an isomorphism
$\mu_*:H_*^{lf}(M)\to H_*^{\Gamma}(X)$ of the homologies defined by
locally finite
chains. The inclusion $ C_*^{lf}(X)^{\Gamma}\subset C_*^{lf}(X)$
defines a
forgetful homomorphism $\nu_*:H^{\Gamma}_*(X)\to H^{lf}_*(X)$.
The composition
$\tau_*=\nu_*\mu_*:H_*^{lf}(X/\Gamma)\to H^{lf}_*(X)$ is called the transfer.
The equivariant cohomology $H^*_{\Gamma}(X)$ of $X$ are defined by the chain complex
$Hom(C_*^{lf}(X)^{\Gamma},\Q)$. There are the dual homomorphisms $\mu^*$, $\nu^*$ and $\tau^*$.

The above definition can be extended to nonregular covering maps. If $p:X\to Y$ is a covering
map then the homomorphism $\mu: C^{lf}_*(Y)\to C^{lf}_*(X)$ can be defined by the formula
$\mu(c)=\sum_{\bar c\in L(c)}\bar c$ where $L(c)$ is the set of all lifts of the singular
simplex $c$. Then the chain complex $\mu(C^{lf}_*(Y))$ defines the
$p$-invariant homology of $X$ and the dual complex defines the $p$-invariant cohomology.
We will use the following notations for them: $H^p_*(X)$ and $H_p^*(X)$.
Note that they are isomorphic to $H^{lf}_*(Y)$ and to $H_c^*(Y)$ respectively by
means of isomorphisms $\mu_*$ and $\mu^*$ induced by $\mu$.

Assuming the diagonal action on $X\times Y$ one
can define an equivariant slant product
$H^n_{\Gamma}(X\times Y)\otimes H_k^{\Gamma}(Y)\to
H^{n-k}_{\Gamma}(X)$
in the same manner as the standard slant product: $[\omega]/[a]=[\omega/a]$ where
$(\omega/a)(b)=\omega(a\times b)$,
$\omega\in Hom(C_n^{lf}(X\times Y)^{\Gamma},\Q)$, $a\in C_k^{lf}(Y)^{\Gamma}$, and
$b\in C^{lf}_{n-k}(X)^{\Gamma}$. One of the features of the equivariant slant product
is the equality $[\omega]/[a]=[\omega']/[a]$ where $[\omega']$ is the image of $[\omega]$ under
the homomorphism $H^*_{\Gamma}(X\times Y)\to H_{\Gamma\times\Gamma}^*(X\times Y)$.
We also recall the main formula for the slant product:
$a\otimes b/c=a\langle b, c\rangle$.

We consider the case of geometrically finite groups $\Gamma$, i.e. groups with
$B\Gamma$ finite complex. Let $X=E\Gamma$ denote the universal cover of $B\Gamma$.
Proper Lipschitz cohomology classes of $\Gamma$ are detected by means of
the following data [CGM]: a finite dimensional complex $P$ with a proper $\Gamma$-action
and a  map $\alpha:X\times P\to\R^N$ satisfying the conditions
\roster
\item{} $\alpha$ is invariant with respect to the diagonal action of $\Gamma$ on
$X\times P$;
\item{} the restriction $\alpha\mid_{x\times P}$ is proper for all $x\in X$.
\endroster
Then the slant product for equivariant (co)homology
$
H^N_{\Gamma}(X\times P)\otimes H^{\Gamma}_i(P)\to H^{N-i}_{\Gamma}(X)
$
defines a homomorphism $$\alpha_{\cap}: H^{\Gamma}_*(P)=H^{lf}_*(P/\Gamma)\to
H^*_{\Gamma}(X)=H^*(B\Gamma)$$ by means of the class
$\alpha^*(\omega)\in H^N_{\Gamma}(X\times P)$
where $\omega\in H^N_c(\R^N)$ is a generator.
Thus, $$\alpha_{\cap}(b)=\alpha^*(\omega)/b.$$

By the definition all classes in $\alpha_{\cap}(H^{\Gamma}_*(P/\Gamma;\Q))\subset
H^*(\Gamma;\Q)$ are called {\it proper Lipschitz cohomology classes} of $\Gamma$
provided the following condition holds

\

(3) the restriction $\alpha\mid_{X\times p}$ is 1-Lipschitz for all $p\in P$.

\

Now we consider the case when $B\Gamma=M$ is a closed orientable $n$-manifold,
$X$
is its universal cover, $\Gamma=\pi_1(M)$.

We note  that if $\Gamma$ acts on $P=X\times\R^k$ in such a way that
the projection $pr:P\to X$ is equivariant then
$H_*^{\Gamma}(P)=H^{lf}_*(X\times_{\Gamma}\R^k)=H_{*-k}(M)$ by means of the Thom
isomorphism for the $\R^k$-bundle $X\times_{\Gamma}\R^k\to M$.

\proclaim{Proposition 1}
Suppose that $\Gamma$ acts on $P=X\times\R^k$ such that the projection
$pr:X\times\R^k\to X$ is equivariant and let $\alpha:X\times P\to\R^{n+k}$ be a map
satisfying the conditions (1)-(3). Then the following are equivalent:
\roster
\item{} $\alpha_{\cap}$ is isomorphism in all dimensions;
\item{} $\alpha_{\cap}\ne 0$;
\item{} $\alpha\mid_{x\times P}$ is essential for all $x\in X$.
\endroster
\endproclaim
\demo{Proof} (1) $\Rightarrow$ (2).
Obvious.

(2) $\Rightarrow$ (3). Since $\alpha$ is $\Gamma$-invariant, it induces a map 
$\hat\alpha:X\times_{\Gamma}P\to\R^{n+k}$, $\alpha=\hat\alpha\circ\bar q$
where $\bar q:X\times P\to X\times_{\Gamma}P$ is the projection to the orbit space. 
The condition (2) implies that 
$\alpha^*(\omega)\ne 0$ and hence $\hat\alpha^*(\omega)\ne 0$ where
$\omega\in H^m_c(R^{n+k})$ is the fundamental class. 
Being the product of a contractible manifold with a Euclidean space, $P$ is homeomorphic
to $\R^{n+k}$.
Since the $P$-bundle $\pi:X\times_{\Gamma}P\to M$ is orientable, the inclusion
$\pi^{-1}(y)\subset X\times_{\Gamma}P$ induces an isomorphism of $(n+k)$-dimensional
cohomology groups with compact supports. Hence the restriction $\hat\alpha$ to $\pi^{-1}(y)$
is essential for all $y\in M$. Therefore the restriction of $\alpha$ to $x\times P$ is
essential for all $x\in X$.

(3) $\Rightarrow$ (1). The essentiality in (3) implies that $\hat\alpha^*(\omega)$ is
a generator of the group $H^{n+k}_c(X\times_{\Gamma}P;\Q)=\Q$. By the proposition 
conditions $P/\Gamma=E$ is the total space of an $\R^k$-bundle $\pi_2:E\to M$. Moreover,
$\pi$ is the pull-back of $\pi_2$ and $\pi_1:X\times_{\Gamma}X\to M$. 
Let $u\in H_c^{n+k}(X\times_{\Gamma}P)$ be the Thom class  of $\pi$. Thus, 
$\hat\alpha^*(\omega)=\lambda u$ for some $\lambda\in\Q\setminus\{0\}$.

First we consider the case when $P=X$. The projection 
$\bar\pi:X\times X\to M\times M$
onto the orbit space of $\Gamma\times\Gamma$-action 
is factored as $\bar\pi=\pi'\circ\bar q$
where $\bar q$ is projection to the 
orbit space of the diagonal action and 
$\pi':X\times_{\Gamma}X\to M\times M$ 
is the induced covering map. Let $N(M)$ denote an
open $r$-neighborhood of $M$ 
imbedded in $X\times_{\Gamma}X$ as the image $\bar q(X)$ of
the diagonal in 
$X\times X$. We will call this imbedding
 a {\it zero section} of $\pi$.
Let $j:N(M)\to X\times_{\Gamma}X$ 
denote the imbedding. If $r$ is small enough, the 
composition $j_1=\pi'\circ j$ 
is an imbedding as well. We may assume that the Thom class
$u$ of $\pi$ is 
supported in $N(M)$. It is known [MS] that $w=j_1^*(u)$ is the Poincare 
dual
in $M\times M$ to the diagonal homology class $M\subset M\times M$. 
Also it is known [MS]
that the slant product with $w$ 
in
$H^n(M\times M)\otimes H_i(M)\to H^{n-i}(M)$ defines the 
Poincare duality for
$M$. The equivariant slant product with 
$\bar w=(\mu^*)^{-1}(w)$ 
in
$$H^n_{\Gamma\times\Gamma}(X\times X)\otimes H_i^{\Gamma}(X)\to H^{n-i}_{\Gamma}(X)$$
also defines an isomorphism. Namely, $\mu^*\circ(\bar w/-)\circ\mu_*=w/-$. The
following
diagram is commutative
$$
\CD
H^*_c(N(M)) @>j^*>> H_c^*(X\times_{\Gamma}X) 
@<\mu^*<< H^*_{\gamma}(X\times X)\\
@Vj_1^*VV  @ V\nu^*VV @ V\nu_1^*VV\\
H^*(M\times M)
 @ <\mu^*<< H_{\pi'}^*(X\times_{\Gamma}X) @<\mu^*<< 
H^*_{\Gamma\times\Gamma}(X\times X)\\
\endCD
$$
We check that $j_1^*=\mu^*\nu^*j^*$. Indeed, 
$j_1^*(a)(c)=a(j_1^{-1}(c|_{j_1(N(M))}))=
q(\mu_*(c)\mid_{N(M)})=a(\mu_*(c))=\mu^*(a)(c)=
\mu^*\nu^*j^*(a)(c)$. Here
$a\in Hom(C^{lf}_*(X\times_{\Gamma}X),\Q)$ with the support
in $N(M)$, $c\in C_*(M)$ with
the mesh $<\epsilon$ where $\epsilon$ is 
such that the $\epsilon$-neighborhood of the
image under $j_1$ of the 
support of $a$ is contained in $j_1(N(M))$. The restriction
$c|_W$ of a chain 
to a neighborhood $W$ is a chain defined by simplices that are contained 
in $W$.
Let $\bar u=(\mu^*)^{-1}(u)$. Then $\nu_1^*(\bar u)=\bar w$.
Therefore 
$\bar u/-:H^{\Gamma}_i(X)\to H^{n-i}_{\Gamma}(X)$ is an isomorphism.
Since 
$\alpha^*(\omega)=\lambda\bar u$, the homomorphism $\alpha_{\cap}$ is an 
isomorphism.

In general case $P=X\times\R^k$ we still can define the diagonal 
imbedding $M\to M\times E$
using a zero section $s:M\to E$ of the bundle $\pi_2$. 
Then similarly one can define the maps
$j$, $j_1$ and $\pi'$. We note that 
the class $w\in H^*(M\times M)$ can be presented as
$\bigoplus_{l}(b_l\otimes c_l)$ 
where $\dim b_l+\dim c_l=n$. Then the image 
under $j_1^*$ of the Thom calss $u$ in 
the general case equals 
$w'=\bigoplus_{l}(b_l\otimes (c_l\cup v))$ where $v$ is the 
Thom class of $\pi_2$.
The slant product with $w'$ in 
$$H^{n+k}_c(M\times E)\otimes 
H_i^{lf}(E)\to H^{n+k-i}(M)$$
defines the composition of the homology Thom isomorphism
$$v\cap -: H^{lf}_i(E)\to H_{i-k}(M)$$ with the Poincare duality. 
Indeed, 
for $a\in H_i^{lf}(E)$ we have
$$w'/a=\bigoplus_{\dim c_l=i-k}b_l\langle c_l\cup v,a\rangle=
\bigoplus_{\dim c_l=i-k}b_l\langle c_l,a\cap v\rangle=w/(a\cap v)=PD(a\cap v).$$
Let $\bar w'=(\mu^*)^{-1}(w')$. 
Then the slant product with $\bar w'$ 
in $$H^{n+k}_{\Gamma}(X\times P)\otimes
H_i^{\Gamma}(P)\to H^{n+k-i}_{\Gamma}(X)$$ is an isomorphism.
The commutative diagram
$$
\CD
H^*_c(N(M)) @>j^*>> H_c^*(X\times_{\Gamma}P) @<\mu^*<< H^*_{\gamma}(X\times P)\\
@Vj_1^*VV  @ V\nu^*VV @ V\nu_1^*VV\\
H_c^*(M\times E) @ <\mu^*<< H_{\pi'}^*(X\times_{\Gamma}P) @<\mu^*<< 
H^*_{\Gamma\times\Gamma}(X\times P)\\
\endCD
$$
imples that the slant product
with $\alpha^*(\omega)=\lambda\bar q^*(u)$
equals the slant product with $\lambda\bar w'$.  Hence $\alpha_{\cap}$ is an isomorphism.\qed
\enddemo

DEFINITION. We say that a manifold $M$ has {\it canonically Lipschitz cohomology}
if there exists 
a map $\alpha$ as in Proposition 1. The map $\alpha$ is called
{\it Lipschitz 
representation} of $H^*(M)$.

Some version of the following was stated in [CGM] (9.A).
\proclaim{Assertion 1} for every closed orientable aspherical manifold $M$ with
the fundamental group $\Gamma$ and the universal cover $X$
there exists a free proper action of $\Gamma$ on $P=X\times\R^k$ for some $k$ with an 
equivariant projection to the first factor and a  
map $\alpha:X\times P\to\R^{n+k}$
satisfying (1)-(2) and such that $\alpha_{\cap}$ is an isomorphism.
\endproclaim
\demo{Proof}
The diagonal $\Gamma$-action defines a bundle
$\pi_1:X\times_{\Gamma}X\to M$ with the fiber $X$ and the structure group $\Gamma$. 
We multiply 
$X\times_{\Gamma}X$ by $\R$
to obtain a topological $\R^{n+1}$-bundle over $M$. 
There is a vector bundle
such that its "Whitney sum" 
with our topological vector bundle gives a trivial
topological vector bundle over $M$. 
It means that there exists a vector bundle 
$\nu:E\to M$ such that the pull-back diagram

$$
\CD
E' @>\pi'>> E\\
@V\nu'VV @V\nu VV\\
X\times_{\Gamma}X @>\pi_1>> M.\\
\endCD
$$
defines a topologically trivial bundle $\pi:E'\to M$, $\pi=\pi_1\circ\nu'$ with
the fiber $\R^{n+k}$.

Let $k=\dim E$. The following pull-back diagram defines $P$ together with an action of 
$\Gamma$ on it
$$
\CD
P @>\bar\nu>> X\\
 @Vq'VV @VqVV\\
E @ >\nu>> M.\\
\endCD
$$
Since $X$ is contractible, $\bar\nu$ is a trivial $\R^k$-bundle, i.e. 
$P$ is homeomorphic to $X\times\R^k$ and $\bar\nu$ is an equivariant projection onto the
first factor.
It is easy to see that the orbit space of the diagonal action of $\Gamma$ on
$X\times P$ is $E'$. Let $p:E'\to\R^{n+k}$ be a trivialization.
We define $\alpha=p\circ\bar q$ where $\bar q:X\times P\to E'$ is the orbit map.
The argument (3) $\Rightarrow$ (1) of Proposition 1 completes the proof.
\qed
\enddemo

To insure condition (3) for the above $\alpha$ one needs the following.

{\it Displacement Bound Condition} [CGM]: 
$
\|p(x)-p(\gamma x)\|\le\|\gamma\|_{S}
$
for all $\gamma\in\Gamma$. 

Here $\|\gamma\|_S$ is the minimal length of $x\in\Gamma$
with respect to a given finite symmetric set of generators $S$.
We assume that 
$S$ is fixed and we let $d_S(x,y)=\|x^{-1}y\|_S$ to denote the word metric on $\Gamma$.

\proclaim{Proposition 2 [CGM]}
Let $q:\Gamma\times P\to P$ be projection to 
the orbit space of the 
diagonal action. Then the following conditions for a map 
$p:P\to\R^N$
are equivalent:
\roster
\item{} $p$ satisfies the Displacement 
Bound Condition;
\item{} The restriction $p\circ q\mid_{\Gamma\times x}$ 
is
1-Lipschitz for all $x\in P$ with respect to the word metric on $\Gamma$.
\endroster

\endproclaim

\demo{Proof} We identify $P$ with $1\times P\subset\Gamma\times P$.
Suppose that $p$  satisfies the Displacement Bound Condition.
Then $\|p \circ q(\gamma_1\times x)-p \circ q(\gamma_2\times x)\|=
\|p \circ q(1\times\gamma_1^{-1}x)-p \circ q(1\times\gamma_2^{-1}x)\|=
\|p (\gamma_1^{-1}x)-p (\gamma_2^{-1}x)\|=
\|p (\gamma_1^{-1}x)-
p((\gamma_2^{-1}\gamma_1)(\gamma_1^{-1}x))\|
\le\|\gamma_2^{-1}\gamma_1\|_{S}=
d_S(\gamma_1,\gamma_2)$.

Now if $p \circ q\mid_{\Gamma\times x}$ is 1-Lipschitz, 
we have
$\|p (x)-p (\gamma x)\|=\|p \circ q(1\times x)-p \circ 
q(1\times\gamma x)\|=
\|p \circ q(1\times x)-p \circ 
q(\gamma^{-1}\times x)\|\le d_S(1,\gamma^{-1})=
\|\gamma\|_S$.
\qed
\enddemo

\proclaim{Proposition 3}
If the above map $p$ is proper and satisfies the Displacement Bound Condition then the map
$p\circ q$ generates a map $\alpha : X\times P\to\R^N$ that satisfies the 
conditions (1)-(3). \endproclaim

\demo{Proof} We may assume that $\Gamma$ is isometrically imbedded in $X$ as an orbit
$\Gamma x_0$ and $\Gamma$ is $K$-dense in $X$. Let $P(\Gamma)$ be
the space of 
measures on $\Gamma$ with finite supports. 
Thus $P(\Gamma)$ consists of formal finite sums $\Sigma\lambda_ix_i$, $\lambda_\ge 0$,
$x_i\in\Gamma$.
Let $1-Lip(\Gamma)$ denote the set
of real-valued 1-Lipschitz functions on $\Gamma$. We endow $P(\Gamma)$ with the 
Kantorovich-Rubinstein metric:
$$d_P(\Sigma\lambda_ix_i,\Sigma\mu_i x_i)=\sup\{|\Sigma(\lambda_i-\mu_i)f(x_i)|
\mid  f\in 1-Lip(\Gamma)\}.$$
Then $\Gamma$ is isometrically imbedded in $P(\Gamma)$ by means of Dirac measures.
Every map $g:\Gamma\to\R^N$ has a natural extension
to 
$\bar g:P(\Gamma)\to\R^N$ by the formula $\bar g(\Sigma\lambda_ix_i)=
\Sigma\lambda_ig(x_i)$. We note that for every 1-Lipschitz map 
$g$ the
extension $\bar g$ is $N$-Lipschitz. Indeed,
$\|\bar 
g(\Sigma\lambda_ix_i)-\bar g(\Sigma\mu_ix_i)\|=
\|\Sigma\lambda_ig(x_i)-
\Sigma\mu_ig(x_i)\|\le$
$$\Sigma_{k=1}^N|\Sigma_i(\lambda_i-\mu_i)g_k(x_i)|\le 
Nd_P(\Sigma\lambda_ix_i,
\Sigma\mu_i x_i)$$ where $g=(g_1,\dots,g_N).$
Thus, $p\circ q$ 
defines a map $\overline{p\circ q}:P(\Gamma)\times P\to\R^N$
such that the restriction to $P(\Gamma)\times x$ is $N$-Lipschitz for all $x$.

A cover of $X$ by $K$-balls $B_K(x)$ centered at points  $x\in\Gamma$ defines a  map 
$\psi :X\to P(\Gamma)$ by the formula
$\psi(x)=\Sigma\lambda_ix_i$ where 
$\lambda_i=d_X(x,X\setminus B_K(x_i))$.
If $d_X(x,y)\le K$, 
then $$d_P(\psi(x),\psi(y))=|
\Sigma(d_X(x,X\setminus B_K(x_i))-
d_X(y,X\setminus B_K(x_i)))f(x_i)|
\le\Sigma d_X(x,y)|f(x_i)|$$
 for some 1-Lipschitz real-valued function $f$ with $f(x)=0$. 
Hence
$|f(x_i)|\le d_X(x_i,x)\le 3K$. Then
$d_P(\psi(x),\psi(y))=
2C(B_K)d_X(x,y)(3K)=Ld_X(x,y)$ where $L=6kC(B_K)$ and $C(B_K)$ is the capacity
of a 
$K$-ball in $\Gamma$. 
Thus, the map $\psi$ is locally $L$-Lipschitz.
Since $X$ is a geodesic metric space, 
$\psi$ is $L$-Lipschitz.

We define $\alpha(x,z)=r\circ 
\overline{p\circ q}\circ(\psi\times 1_P)$
where $r:\R^N\to\R^N$ is a 
rescaling that makes $\alpha$ 1-Lipschitz on $X\times x$, $x\in P$.
\qed
\enddemo

A variation of a continuous function $f:X\to\R$ on a metric space $X$
at point $x\in X$ in an $R$-ball is
defined as follows
$$
var_R(f)(x)=\max\{|f(x)-f(x')|\ \mid x'\in B_R(x)\}.
$$
We recall that the {\it Higson corona} $\nu X$ of a proper metric space
$X$ is the corona $\bar X\setminus X$ of the compactification
of $X$ defined by the ring of bounded functions with the property
$
\lim_{dist(x,x_0)\to\infty}var_R(f)(x)=0$
for every $R$ and fixed (every) $x_0\in X$.

Suppose that a finitely generated group $\Gamma$ acts properly on a locally
compact space $X$. We define a $\Gamma$-variation of a continuous function $f:X\to\R$
at point $x\in X$ in an $R$-ball as follows
$$
\Gamma-var_R(f)(x)=\max\{|f(x)-f(x')|\ \mid\ \ x'=\gamma x,\ \|\gamma\|\le R\}.
$$
The algebra of bounded continuous functions $f$ on $X$ with the property that 
for every $R>0$ the
$\Gamma$-variation $\Gamma-var_R(f)(x)$ tends to zero as $x$ approaches infinity
defines a compactification $\tilde X$ of $X$ which we will call the 
{\it Hurder compactification} (cf. [Hu]). The corona 
$\partial_{\Gamma}X=\tilde X\setminus X$
will be called the {\it Hurder corona}.

Let $\overline{\Gamma\times\R}$ be the 
compactification of 
$\Gamma\times\R$ by the suspension $\Sigma(\nu\Gamma)$ of
the Higson corona.
When $X/\Gamma$ is compact the space 
$\overline{\Gamma\times\R}\times_{\Gamma}X$ 
is a compactification of the space $(\Gamma\times\R)\times_{\Gamma}X=X\times\R$
with the corona  $\Sigma(\nu\Gamma)\times_{\Gamma}X$.

\proclaim{Proposition 4} Suppose that the action of a finitely generated group $\Gamma$ 
on a space
$X$ is free. Then
\roster
\item{} The Hurder compactification of $X$ is the maximal
$\Gamma$-equivariant compactification for which the corona is the fixed point set;
\item{} $\partial_{\Gamma}X=\nu\Gamma\times_{\Gamma}X$, provided the action is cocompact;
\item{} There is a map of coronas $\partial_{\Gamma}(X\times\R)\to
\Sigma(\nu\Gamma)\times_{\Gamma}X$ which extends continuously to $X\times\R$ 
by the identity map.
\endroster
\endproclaim

\demo{Proof}

(1) Clearly the action of $\Gamma$ on the Hurder corona is trivial.
It remains to show that every bounded function $f:X\to\R$ satisfying the property
$f-f\circ\gamma\in C_0(X)$ for all $\gamma\in\Gamma$ belongs to the Hurder algebra,
i.e. $\lim_{x\to\infty}\Gamma-var_R(f)(x)=0$. The later follows from the fact that
$\lim_{x\to\infty}(f(x)-f(\gamma x))=0$ and the faintness of the $R$-ball $B_R(e)$ in
$\Gamma$.

(2) The proof is contained in [Hu], Lemma 3.1. We add one comment here. This corona is obtained
by taking the Higson corona of every orbit $\Gamma x$, $x\in X$, supplied by the right invariant
word metric $\rho$. Then it fibered over $M$ with the fiber $\nu\Gamma$ as it stated in (2). 
We note the right invariant metric on $\Gamma$ is defined as
$d^R(\gamma',\gamma)=\|\gamma'\gamma^{-1}\|_S$ and for every orbit $\Gamma x$ the formula
$\rho(\gamma'x,\gamma x)=d^R(\gamma',\gamma)$ defines a metric on $\Gamma x$ 
together with an isometry to $(\Gamma,d^R)$.
As a metric space $(\Gamma,d^R)$ is isometric to $\Gamma$ with the left invariant metric 
by taking the inverse of the elements. Hence
the Higson corona of every orbit $\Gamma x$ is just $\nu\Gamma=\nu X$. We note that the left action
of every $\gamma\in\Gamma$ defines a translation on $\Gamma x$ bounded by $\|\gamma\|_S$ for every
$x\in X$. This makes the action of $\Gamma$ on the corona 
$\nu(\Gamma)\times_{\Gamma}X$ trivial.

(3) We consider a compactification of $X\times\R$ by taking the Higson corona of
$(\Gamma x,d^R)\times\R$ for every orbit $\Gamma x$. 
It is easy to check that the corona of this compactification is $\nu(\Gamma\times\R)\times_{\Gamma}X$.
The action of $\Gamma$ on this corona
is trivial by the reason given in (2). By (1) there is a domination
$\beta:\partial_{\Gamma}(X\times\R)\to\nu(\Gamma\times\R)\times_{\Gamma}X$.
The obviuos domination $\nu(\Gamma\times\Z)\to\Sigma(\nu\Gamma)$ defines a map
$\nu(\Gamma\times\Z)\times_{\Gamma}X\to\Sigma(\nu\Gamma)\times_{\Gamma}X$
which together with $\beta$  gives us a required map.\qed
\enddemo

\proclaim{Proposition 5}
Let $\phi:\partial_{\Gamma}X\to S^{n-1}$ 
be an essential map of the Hurder corona of a contractible 
$n$-manifold $X$ supplied with a 
$\Gamma$-action such that
$\delta([\phi])\ne 0$
, 
$\delta:\check{H}^{n-1}(\partial_{\Gamma}X)\to H^n_c(X)$. 
Then there exists an essential map $p:X\to\R^n$ satisfying the 
Displacement Bound Condition.
\endproclaim

\demo{Proof}
There is an extension $\bar\phi:W\to S^{n-1}$ of $\phi$ to an open neighborhood of
$\partial_{\Gamma}X\subset W\subset\tilde X$. Let $r_0>0$ be such that 
$X\setminus W\subset B_{r_0}(x_0)$ for some $x_0\in X$. Let $S=S^{-1}$ be
a finite generating set of $\Gamma$. We define
$$\nu(t)=\max\{\rho(\bar\phi(x),\bar\phi(sx))\ \mid\ s\in S,\ \|x\|\ge t\}$$
where $\rho$ is the standard metric on the unit sphere and $\|x\|=d_X(x,x_0)$.
Since $\lim_{\|x\|\to\infty}\Gamma-var_1(\bar\phi)(x)=0$, we have
$\lim_{t\to\infty}\nu(t)=0$. There exists a smooth function $f:\R_+\to\R_+$ such that
$f(t)=0$ for $t\le r_0$, $f(t)\le1/\nu(t)$ and $f'(t)<1$ for all $t$. We define
$p:X\to\R^n$ in polar coordinates as $p(x)=(1/2f(\|x\|),\bar\phi(x))$.
Note that $\|x\|+1\ge\|sx\|\ge\|x\|-1$ for every $s\in S$. Then
$\|p(x)-p(sx)\|\le$  

$$1/2|f(\|sx\|)-f(\|x\|)|+
1/2f(\|x\|)\sin(\rho(\bar\phi(x),\bar\phi(sx))
\le 1/2(1+f(\|x\|)\nu(\|x\|))\le 1.$$ If $\|\gamma\|_S=k$ and 
$\gamma=s_1\dots s_k$, we
obtain $$\|p(x)-p(\gamma x)\|\le
\|p(x)-p(s_kx)\|+\dots+\|p(s_2\dots p_kx)-p(s_1\dots p_k x)\|\le k=\|\gamma\|_S.$$
Thus the Displacement Bound Condition is satisfied.

To show that $p$ is essential we consider a map $g:X\to B^n$ defined as
$g(x)=p(x)/(1+p(x))$. The map $g$ admits a continuous extension $\bar g$ to the
Hurder corona by means of $\phi$. Then the diagram
$$
\CD
\check H^{n-1}(\partial_{\Gamma}X) @>\delta>> \check H^n(\tilde X,
\partial_{\Gamma}X)\\
@ A\phi^*AA @ A\bar g^*AA\\
H^{n-1}(S^n) @>\delta>> H^n(B^n,S^{n-1})\\
\endCD
$$
implies that $g$ is essential. Therefore $p$ is essential.\qed
\enddemo

We recall that $\beta Y$ is the standard notation for the
Stone-\v Cech compactification of $Y$. We denote by $\check\beta Y=
\beta Y\setminus Y$ the Stone-\v Cech corona.

\proclaim{Theorem 1}
The following conditions for an $n$-dimensional closed orientable aspherical manifold
$M$ with the fundamental group $\Gamma=\pi_1(M)$ and the universal
cover $X$ are equivalent:
\roster
\item{} There exist $k$ and a $\Gamma$-action on $X\times\R^k$ with
an equivariant projection onto $X$ such that the corona of the
Stone-\v Cech compactification of $X\times\R^k$ admits a $\Gamma$-invariant
essential map onto $n+k-1$-dimensional sphere
$f:\check\beta(X\times\R^k)\to S^{n+k-1}$.
\item{} There exist $k$ and a $\Gamma$-action on $X\times\R^k$ with
an
equivariant projection onto $X$ such that the boundary homomorphism
$\delta: \check H^{n+k-1}(\partial_{\Gamma}(X\times\R^k);\Q)
\to H_c^{n+k}(X\times\R^k;\Q)$
for the Hurder corona is an epimorphism;
\item{} There exist a $\Gamma$-action on $X\times\R^k$ for some $k$ with
an equivariant projection onto $X$ and a proper essential map
$p:X\times\R^k\to\R^{n+k}$
satisfying the Displacement Bound Condition;
\item{} The manifold $M$ has the rational cohomology canonically Lipschitz.
\endroster
\endproclaim
\demo{Proof}
(1) $\Leftrightarrow$ (2). Let $f$ be as in (1), then it defines a
compactification of $X\times\R^k$ with the corona $S^{n+k-1}$ and the trivial
$\Gamma$-action on it. In view of Proposition 4 (1) this compactification is
dominated by the Hurder compactification. Hence we have the following
commutative diagram:
$$
\CD
\check H^{n+k-1}(\check\beta(X\times\R^k)) @>\delta_1>> H^{n+k}_c(X\times\R^k)\\
@ Af_1^*AA  @ AidAA\\
\check H^{n+k-1}(\partial_{\Gamma}(X\times\R^k)) @ >\delta>> H_c^{n+k}(X\times\R^k)\\
@ Af_1^*AA  @ AidAA\\
\check H^{n+k-1}(S^{n+k-1}) @ >\delta_2>> H_c^{n+k}(X\times\R^k)\\
\endCD
$$
By the Calder-Siegel theorem [CS] $\delta_1$ is an isomorphism for $n+k>2$.
Since the map $f=f_2\circ f_1$ is essential, $\delta_1\circ f^*$ is a 
rational isomorphism. Hence $\delta$ is an epimorphism.

Assume that the condition (2) is satisfied. We may assume that $n+k-1$ is odd.
Then an Eilenberg-Maclane complex $K(\Q,n+k-1)$ may be regarded as an infinite
telescope of $(n+k-1)$-spheres.
Then there exists a map $f_2:\partial_{\Gamma}(X\times\R^k)\to S^{n+k-1}$
such that $\delta\circ f^*_2\ne 0$. Let $f_1:\check\beta(X\times\R^k)\to
\partial_{\Gamma}(X\times\R^k)$ be a map that comes from the universality of the
Stone-\v Cech compactification. Then we have the above diagram which 
together with Calder-Siegel theorem imply that $f=f_2\circ f_1$ is essential.

(2) $\Leftrightarrow$ (3).
 If (2) holds for $k$ then it holds for all $k'\ge k$.
We take $k$ such that $n+k-1$ is odd.  Then like in the argument for (2) $\Leftrightarrow$ (3)
there exists a map $\phi:\partial_{\Gamma}(X\times\R^l)\to S^{n+k-1}$
such that $\delta([\phi])\ne 0$. Then (3) holds by Proposition 5.

Let $p:X\times\R^k\to\R^{n+k}$ be as in (3). 
We define a map 
$f:X\times\R^k\to B^{n+k}$ as $f(x)=p(x)/(1+\|p(x)\|)$. It is easy 
to show that
$f$ satisfies the condition $\lim_{x\to\infty}\Gamma-var_R(f)(x)=0$. 
Then $f$ is extendible to the Hurder corona. 
Then the essentiality of the map 
$f:(\widetilde{X\times\R^k},\partial_{\Gamma}(X\times\R^k))\to 
(B^{n+k},\partial B^{n+k})$ and the cohomology exact sequence of pair imply (2).

(3) $\Leftrightarrow$ (4).
By Proposition 3 there exists
$\alpha:X\times 
P\to\R^{n+k}$ with $P=X\times\R^k$ satisfying (1)-(3).
Since $p$ is essential, 
in view of Proposition 1 $\alpha$ is a Lipschitz representation of $H^*(M)$.

If $\alpha:X\times P\to\R^{n+k}$ gives a Lipschitz 
representation of $H^*(M)$,
the restriction $\alpha'=\alpha|_{\Gamma x_0\times P}$ 
is factored through
a map $p:\Gamma\times_{\Gamma}P=P\to\R^{n+k}$, 
$\alpha'=q\circ p$, which 
satisfies the Displacement Bound Condition 
by Proposition 2.\qed
\enddemo

DEFINITION. {\it We say that the universal cover $X$ of a closed aspherical manifold $M$
with the metric induced from $M$
satisfies the condition CP if it admits a rationally acyclic $\pi_1(M)$-equivariant
Higson dominated compactification } [CP], [D].

\proclaim{Theorem 2}
Suppose that the universal cover $X$ of a closed orientable aspherical $n$-manifold $M$
satisfies the conditions CP. Then the cohomology groups $H^*(M;\Q)$ are canonically
Lipschitz.
\endproclaim
\demo{Proof}
Let $Y$ be the above compactification of $X$ and let $Z=Y\setminus X$.
Since $H^*(Y;\Q)=0$, we have $H^*(Z;\Q)=H^*(S^{n-1};\Q)$ for the reduced \v{C}ech cohomology.
Moreover, since the action of $\Gamma=\pi_1(M)$ is trivial on $H^*(X;\Q)$ and
$\delta:H^*(Z;\Q)\to H^{*+1}(X;\Q)$ is an equivariant isomorphism, the action of
$\Gamma$ on $H^*(Z;\Q)$ is trivial.

Let $\Sigma Z$ be the suspension. We consider a locally trivial fibration 
$\pi:\Sigma Z\times_{\Gamma}X=S\to M$
with the fiber $\Sigma Z$ induced by the projection $pr:\Sigma Z\times X\to X$.
The term $E_2^{0,n}$ in the spectral sequence of this fibration (with local coefficients)
is isomorphic to $H^n(\Sigma Z;\Q)=\Q$ since the action of $\Gamma$ is trivial on
$H^n(\Sigma Z;\Q)$. By dimensional reason the term $E_2^{0,n}$ is not transgressive.
Therefore the inclusion of the fiber $\Sigma Z\subset S$
induces an epimorphism of $n$-dimensional rational cohomology groups.
 
The space $S$ is contained in the
space $\bar E=\widetilde{X\times\R}\times_{\Gamma}X$ where $\widetilde{X\times\R}$
is the natural compactification of $X\times\R$ by $Z\ast S^0=\Sigma Z$.
The exact sequence of pair and the inclusion
$(\widetilde{X\times\R},\Sigma Z)\subset (\widetilde{X\times\R}\times_{\Gamma}X,S)$ defines the
commutative diagram for rational cohomology:
$$
\CD
\check H^n(S) @>\bar\delta>> H^{n+1}_c((X\times\R)\times_{\Gamma}X)\\
@ Vi^*VV @ V j^*VV\\
\check H^n(\Sigma Z) @ >\delta_1>> H_c^{n+1}(X\times\R).\\
\endCD
$$
We proved that there is $a\in\check H^n(S)$ such that $i^*(a)$ is
a generator of $\check H^n(\Sigma Z)=\Q$. Then
$\bar\delta(a)$ is a generator of $H^{n+1}_c((X\times\R)\times_{\Gamma}X)=\Q$.
The later equality follows from the fact that $(X\times\R)\times_{\Gamma}X$ is the total
space of a topological orientable $n+1$-dimensional vector bundle over $M$.
Then $\bar\delta(a)$ is the rational Thom class for two $X\times\R$-bundles over $M$ generated by
the projections of $X\times\R\times X$ onto both factors $X$.
Therefore the restriction of $\bar\delta(a)$ onto any fiber of the first bundle is nontrivial.
Let $x_0\in M$, then the fiber for the first bundle at $x_0$ can be presented as
$X\times\R=\Gamma x_0\times_{\Gamma}(\R\times X)\subset\widetilde{X\times\R}
\times_{\Gamma}X$.

We compactify $X\times\R\times_{\Gamma}X$ by $\Sigma\nu X\times_{\Gamma}X$.
This also gives a compactification with the above corona of 
a fiber from the first bundle $X\times\R=\Gamma x_0\times_{\Gamma}\R\times X$.
It is easy to check that the action of $\Gamma$ on $X\times\R$ give rise to a trivial action
on the corona. In view of proposition 4 the above compactification is Harder dominated.

Since $Z$ is Higson dominated, there is a map of the Higson corona
$\xi:\nu X=\nu\Gamma\to Z$ that extends to $X$
continuously by the identity map. This map generates a map
$\xi:\Sigma\nu\Gamma\times_{\Gamma}X\to\Sigma Z\times_{\Gamma}X=S$ fixed on $X\times\R$.
By Proposition 4 there is a map
$\xi_1:\partial_{\Gamma}(X\times\R)\to\Sigma\nu\Gamma\times_{\Gamma}X$ 
that extends continuously to $X\times\R$ by
the identity map. The following diagram
$$
\CD
\check H^n(S) @>\delta'>> H^{n+1}_c(X\times\R)\\
@ V(\xi\xi_1)^*VV @ VidVV\\
\check H^n(\partial_{\Gamma}(X\times\R)) @ >\delta>> H_c^{n+1}(X\times\R).\\
\endCD
$$
and the fact that $\delta'$ is an epimorphism implies that
$\delta$ is an epimorphism and hence the condition (2)
of Theorem 1 is satisfied.
\qed
\enddemo

\proclaim{Lemma 1}
Let $X$ be the universal cover $X$ of a closed aspherical
Riemannian $n$-dimensional manifold $M$ whose rational cohomologies are 
canonically Lipschitz and let
$f:X\times P\to\R^N$, $P=X\times\R^k$, be their Lipschitz presentation. Then
there exists a function $c(r)$ such that
$\diam(f^{-1}(B_r(0))\cap X\times
v)\le c(r)$ for all $v\in P$
where 
$B_r(0)\subset\R^N$ is the $r$-ball centered at the origin.
\endproclaim

\demo{Proof} 
Let $\bar q:X\times P\to X\times_{\Gamma}P$ be the orbit map, then $f=\hat f\circ\bar q$
where $\hat f$ is a proper map. Hence $Y=\hat f^{-1}(B_r(0))$ is compact.
Let $\sU$ be a finite cover of $M$ by evenly covered open maps with respect to the
universal cover $q:X\to M$. We can subdivide $Y$ in a finite union 
$Y=\cup_{U\in\sU}Y_U$ of
compact sets such that $Y_U\subset\pi^{-1}(U)$ where
$\pi:X\times_{\Gamma}P\to M$ is a fibration induced by the projection of $X\times P$ onto
the first factor. Let $\pi_2:X\times_{\Gamma}X\to M$ be a fibration induced by the projection
$pr_2:X\times X\to X$ onto the second factor and
let $h_U:\pi_2^{-1}(U)\to X$ be a trivialization defined by means of a lift $s_U$ of $U$ to
$X$. Let $pr':X\times_{\Gamma}P\to X\times_{\Gamma}X$ be a map induced by the
projection $pr:P\to X$.
Let $d'$ be a metric on $M$, We consider the product metric
$\bar d$ on 
$U\times X$ and define
$$
c(r)=\sum_{U\in\sU}\diam_{\bar d}(h_U\times 1_U)(pr'Y_U).
$$

Then 
$$\diam_{d_X}f^{-1}(B_r(0))\cap(X\times v)=
\diam_{d_X}Y\cap\pi^{-1}(q(pr(v)))
\le$$

$$\diam_{d_X}pr'(Y)\cap\pi_2^{-1}(q(pr(v)))
\le$$

$$\sum_{U\in\sU}\diam_{d_X}pr'(Y_U)\cap\pi_2^{-1}(q(pr(v)))\le
\sum_{U\in\sU}
\diam_{\bar d}(h_U\times 1_U)(pr'(Y_U))=c(r).$$
\qed\enddemo

\head \S 3 Expanders\endhead

Let $X$ be a finite graph, we denote by $V$ the set of vertices and by $E$
the set of edges in $X$. We will identify the graph $X$ with its set of vertices
$V$. Every graph is a metric space with respect to the natural metric where every
edge has the length one. For a subset $A\subset X$ we define the boundary
$\partial A=\{x\in X\mid dist(x,A)=1\}$. Let $|A|$ denote the
cardinality of $A$.

DEFINITION [Lu]. An {\it expander} with a conductance number $c$ and the degree
$d$ is an infinite sequence of finite graphs $\{X_n\}$ with the degree $d$ such that
$|X_n|$ tends to infinity and for every $A\subset X_n$ with $|A|\le |X_n|/2$
there is the inequality $|\partial A|\ge c|A|$.

Let $X$ be a finite graph, we denote by $P$ all nonordered pairs of distinct
points in $X$. For every nonconstant map $f:X\to l_2$ to the Hilbert space
we introduce the number
$$
D_f=
\frac{\frac{1}{|P|}\sum_{\{x,y\}\in P}\|f(x)-f(y)\|^2}{\frac{1}{|E|}
\sum_{\{x,y\}\in E}\|f(x)-f(y)\|^2}.
$$
If $X$ is a graph with the degree $d$ and with $|X|=n$, then $|P|=n(n-1)/2$
and $|E|=dn/2$.
The following lemma can be derived from [Ma, Proposition 3].
It also can be obtain from the equality
$$
\lambda_1(X)=\inf\{\frac{\|d f\|^2}{\|f\|^2}\mid \sum f(x)=0\}
$$
for the first positive eigenvalue of the Laplacian on $X$ and the Cheeger's
inequality (see Proposition 4.2.3 in [Lu]).

\proclaim{Lemma 2} Let $\{X_n\}$ be an expander. Then there is a constant
$c_0$ such that $D_{f_n}\le c_0$ for all $n$ for all possible maps
 $f_n:X_n\to l_2$  to the Hilbert space $l_2$.
\endproclaim
\proclaim{Corollary 1}
For every sequence of 1-Lipschitz maps $f_n:X_n\to l_2$ there is the inequality
$\frac{1}{|P_n|}\sum_{P_n}\|f_n(x)-f_n(y)\|^2\le c_0$ for every $n$.
\endproclaim
\demo{Proof}
In the case of 1-Lipschitz map we have
$${\frac{1}{|E|}\sum_{\{x,y\}\in E}\|f(x)-f(y)\|^2}\le 1.$$

Then the required inequality follows.\qed
\enddemo
\proclaim{Corollary 2}
Assume that for a sequence of 
1-Lipschitz maps $f_n:X_n\to l_2$
we have $\sum_{x\in X_n}f_n(x)=0$ for every $n$, 
then
$$
\frac{1}{|X_n|}\sum_{x\in X_n}\|f_n(x)\|^2\le c_0/2
$$
for all $n$.

\endproclaim
\demo{Proof}
We assume that the sequence $\{X_n\}$ is 
enumerated by a subsequence
on $\N$ such that $|X_n|=n$.

According to Corollary 1 we have

$$
n(n-1)c_0/2=|P_n|c_0\ge\sum_{P_n}\|f_n(x)-f_n(y)\|^2_2=
\sum_{P_n}\|f_n(x)-f_n(y)\|^2_2+
(\sum_{X_n}f(x_n))^2=$$

$$
\sum_{P_n}(f_n(x)-f_n(y))^2+\sum_{X_n}f_n(x)^2+\sum_{P_n}2(f_n(x),f_n(y))
=
|X_n|\sum_{X_n}f_n(x)^2=n\sum_{X_n}\|f_n(x)\|^2_2.$$

Then the required inequality follows.\qed
\enddemo

\proclaim{Corolary 3}
Assume that for a sequence of 1-Lipschitz maps 
$f_n:X_n\to l_2$
we have $\sum_{x\in X_n}f_n(x)=0$ for every $n$, 
then there is $R$
such that $|f_n^{-1}(B_R(0))|>|X_n|/2$ for all $n$.
\endproclaim

\demo{Proof}
Take $R>\sqrt{c_0}$. Assume that $|f_n^{-1}(B_R(0))|\le |X_n|/2$

By Corollary 2 we have the contradiction:

$$\frac{c_0}{2}|X_n|\ge\sum_{x\in X_n}\|f(x)\|^2\ge R^2|X_n|/2>\frac{c_0}{2}|X_n|.
$$
\qed
\enddemo

Let $\{X_n\}$ be an expander and let $X$ be a metric space. 
We call a
sequence of 1-Lipschitz maps $\xi_n:X_n\to X$ a 
{\it coarse quasi-embedding}
of an expander $\{X_n\}$ 
into $X$ if there is a positive function
$\rho:\R^2_+\to\R_+$ such that 
$$\lim_{t\to\infty}\frac{\rho(r,t)}{t}=0$$ 
for every $r$, and

$$
|\xi_n^{-1}(B_r(x))|<\rho(r,|X_n|)
$$
for all $n\in\N$ and for all $x\in X$.

When $X$ is a space of bounded geometry, for instance $X$ is a finitely generated group
with the word metric, then this condition is equivalent to the following 
condition defined in [HLS]
$$
\lim_{n\to\infty}\frac{|\xi_n^{-1}(x)|}{|X_n|}=0.
$$

\proclaim{Lemma 3}
Suppose a metric space $(X,d)$ has a bounded geometry. 
Let $f:X\times\R^N\to\R^N$
be a map satisfying the conditions:

\roster
\item{} the restriction $f\mid_{X\times w}:X\times w\to\R^N$ 

is 1-Lipschitz for all $w\in\R^N$;
\item{} the restriction 
$f\mid_{x\times \R^N}:x\times\R^N\to\R^N$ is essential
for all 
$x$;
\item{} there is a function $c(r)$ such that 

$\diam_{d}(f^{-1}(B_r(0))\cap(X\times w))\le c(r)$ for all $w\in\R^N$.

\endroster
Then no expander is coarsely quasi-embeddable in $X$.
\endproclaim

\demo{Proof}

Assume that $\{\xi_n:X_n\to X\}$ is a coarse quasi-embedding of an expander

$\{X_n\}$ into $X$ and let $\rho$ be a corresponding function.

Given $n$ we show that there is $w_n\in\R^N$ such that

$$
\sum_{v\in X_n}f(\xi_n(v),w_n)=0.
$$
We consider a map 
$F:\R^N\to\R^N$ defined as 
$F(y)=\frac{1}{|X_n|}\sum_{v\in X_n}f(\xi_n(v),y)$. 
The condition (1) implies
$\|f(\xi_n(v),y)-f(x_0,y)\|\le d(\xi_n(v),x_0)$. 
Hence
$$
\|F(y)-f(x_0,y)\|=
\frac{1}{|X_n|}\|\sum_{v\in X_n}(f(\xi_n(v),y)-f(x_0,y))\|\le
$$
$$\le\frac{1}{|X_n|}
\sum_{v\in X_n}\|(f(\xi_n(v),y)-f(x_0,y))\|\le 
d(\xi(v_n),x_0),
$$
where 
$d(\xi(v_n),x_0)=\max\{d(\xi(v),x_0)\mid v\in X_n\}$.
Thus the maps 
$F$ and $f\mid_{x_0\times\R^N}$ are in a finite distance.

Therefore they are properly homotopic. Then by the condition (2)
the 
map $F$ is essential. Hence there is $w_n$ with $F(w_n)=0$.
Let 
$\epsilon_n:X\to X\times w_n$ be a natural identification.

We apply Corollary 3 to the sequence of maps
$$(f\mid_{X\times w_n})\circ\epsilon_n\circ\xi_n:X_n\to\R^N$$

to obtain a number $R>0$ with the property

$$
|((f\mid_{X\times w_n})\circ\epsilon_n\circ\xi_n)^{-1}(B_R(0))|\ge |X_n|/2
$$

for all $n$. By the condition (3) 
the preimage 
$(f\mid_{X\times w_n})^{-1}(B_R(0))$ is 
contained in a ball $B_{c(R)}(x_n)$
in $X\times w_n$ of radius 
$c(R)$ centered at some point $x_n$.

Then $$|\xi_n^{-1}\epsilon_n^{-1}(B_{c(R)}(x_n))|\ge
|X_n|/2.$$

Thus,
$$
|X_n|/2\le |\xi_n^{-1}(B_{c(R)}(x_n))|<\rho(c(R),|X_n|)
$$
for all $n$. 
This contradicts with the condition
$\lim_{t\to\infty}\frac{\rho(c(R),t)}{t}=0$.
\qed
\enddemo
Proposition 1, Lemma 1 and Lemma 3 imply the following.

\proclaim{Theorem 3}
Suppose that a closed orientable aspherical manifold $M$
has canonically Lipschitz cohomology. Then no expander admits a coarse quasi-embedding
into its fundamental group. 

\endproclaim

\enddocument